\documentclass{article}
\usepackage[utf8]{inputenc}
\usepackage[letterpaper, total={6in, 8.5in}]{geometry}
\usepackage{graphicx}
\graphicspath{ {./images/} }
\usepackage{wrapfig}
\usepackage{enumitem}
\usepackage{amsmath}
\usepackage{amssymb}
\usepackage{marginnote}
\usepackage{comment}
\usepackage{hyperref}
\usepackage{url}
\usepackage{esint}
\usepackage{multirow}
\usepackage{comment}
\usepackage{caption}
\usepackage{subcaption}
\newcommand{\transpose}{^\intercal}
\usepackage{mathtools}

\title{RESEARCH REPORT---Summer 2023  \\Describing Chaotic Systems\footnote{\begin{center}DISTRIBUTION STATEMENT A.
Approved for public release:
distribution unlimited.\\[0.25cm]

Note: The derivations in this paper were converted into Appendices B.1, B.2, and B.4 in \cite{le2024exploring}.\end{center}}}
\author{\it{Brandon Le}, Naval Research Laboratory, Code 5580, NRL SEAP Program\\
Mentor, Dr. Ira S. Moskowitz}
\date{}

\begin{document}

\maketitle
\vspace{-1.0cm}
\tableofcontents
\listoffigures

\vspace{0.5cm}
\begin{abstract}
    In this paper, we discuss the Lyapunov exponent definition of chaos and how it can be used to quantify the chaotic behavior of a system. We derive a way to practically calculate the Lyapunov exponent of a one-dimensional system and use it to analyze chaotic behavior of the logistic map, comparing the $r$-varying Lyapunov exponent to the map's bifurcation diagram. Then, we generalize the idea of the Lyapunov exponent to an $n$-dimensional system and explore the mathematical background behind the analytic calculation of the Lyapunov spectrum. We also outline a method to numerically calculate the maximal Lyapunov exponent using the periodic renormalization of a perturbation vector and a method to numerically calculate the entire Lyapunov spectrum using QR factorization. Finally, we apply both these methods to calculate the Lyapunov exponents of the H\'enon map, a multi-dimensional chaotic system. 
\end{abstract}

\newpage

\section{Introduction}

There are many definitions of chaos, but the one defined by Lorenz, the father of chaos theory, is that a system is chaotic if it is ``sensitively dependent on initial conditions''. According to Lorenz, initial conditions are not necessarily the conditions of a system when it is created; rather, they can be the conditions at the beginning of any stretch of time that is to be investigated. Furthermore, sensitive dependence implies more than simply an increase in the difference of two states over time. For example, a system where an initial difference of one unit eventually increases to a hundred units and where a difference of a hundredth of a unit eventually increases to one unit is not considered chaotic \cite{lorenz}.

One simple model that can be chaotic is the logistic map. The logistic map is a one-dimensional discrete-time dynamical system with the quadratic recurrence equation
\begin{equation}
    x_{n+1} = rx_n\left(1-x_n\right) \label{eq:1}
\end{equation}
where $r$ is a constant parameter. This map is capable of ``very complicated behavior'' despite its simplistic definition \cite{logistic}. We can model the logistic map in R using the following code:
\begin{verbatim}
    r <- 4
    n <- 10000
    
    k <- numeric(n)
    x <- numeric(n)
    x_0 <- 0.1
    
    k[1] <- 1
    x[1] <- r * x_0 * (1 - x_0)
    
    for (i in 1:(n-1)) {
      k[i+1] = i+1
      x[i+1] <- r * x[i] * (1 - x[i])
    }
    
    plot(k, x, type = "p", pch = 16, 
        xlab = "k", ylab = "x", 
        cex = 0.25, main = "Logistic Map")
\end{verbatim}
Here, we set the parameter $r$ to 4 (which we will later see results in chaotic behavior) and $n$, the number of steps, to 10000. The vector $k$ represents all the steps, and $x$ is a vector of all $x_k$. Because R indexes vectors from 1, we set the initial condition $x_0$ and use Equation \ref{eq:1} to define $x_1$. Then, we loop through all the steps and apply Equation $\ref{eq:1}$ to calculate all the values in $x$. Finally, we plot $x$ vs. $k$ as a scatter plot. Although we can visually see varying complexities in the plot of the logistic model for different values of $r$ and see varying sensitivity by tweaking the initial conditions, it will be useful to have a more rigorous definition of chaos that can quantify how chaotic a given system is.

\section{Quantifying Chaos}

This more specific definition of chaos builds on Lorenz's idea, chaos being defined by Alligood et al. as ``a Lyapunov exponent greater than zero'' \cite{alligood}. In this regard, the Lyapunov exponent can be considered a mathematical indication of chaos. The Lyapunov exponent describes the rate of separation of two trajectories initially separated by some small value $\varepsilon$. For example, consider a one-dimensional discrete-time dynamical system that takes in an input $x_0$ and outputs $x_1 = f(x_0)$. Then, $x_n = f^n(x_0)$. Using the Lyapunov exponent $\lambda$, the following equation is satisfied:
\begin{equation}
    |f^n(x_0 + \varepsilon) - f^n(x_0)| \approx \varepsilon e^{\lambda n} \label{eq:2}
\end{equation}
where we have assumed that the distance between two initially close states grows exponentially with the number of steps $n$. Based on this definition, if $\lambda$ is greater than zero, then any difference in the initial conditions of two systems will result in larger and larger deviations in their trajectories as time goes on. In other words, no matter how small the difference in initial conditions is, the difference in the systems' trajectories will eventually become substantial due to the nature of the exponential function. This demonstrates that a positive $\lambda$ implies a chaotic system by Lorenz's definition. Similarly, a negative $\lambda$ implies that any difference in initial conditions will eventually go to 0 and the two trajectories will converge. This indicates that a negative $\lambda$ is indicative of a dissipative system or an attraction to a stable point or orbit. If $\lambda=0$, any difference in initial conditions will be maintained, neither being magnified nor reduced, meaning the system is in a steady state. The more positive $\lambda$ is, the more quickly a difference in initial conditions will be magnified and the more chaotic a system is. The more negative $\lambda$ is, the more quickly a difference in initial conditions will become negligible and the more dissipative a system is. Thus, the value of $\lambda$ reflects how quickly the system will become stable or unpredictable and can be considered a measurement of chaos, beyond just an indication \cite{wolf}. 

Solving for $\lambda$ from Equation $\ref{eq:2}$,
\begin{equation}
    \lambda\approx\frac{1}{n}\ln\left(\frac{|f^n(x_0 + \varepsilon) - f^n(x_0)|}{\varepsilon}\right) \label{eq:3}
\end{equation}
The ideal Lyapunov exponent that most accurately measures chaos is measured as the number of steps $n$ goes to infinity and the distance between the initial states $\varepsilon$ is infinitesimal:
\begin{align*}
    \lambda &= \lim_{n\to\infty}\lim_{\varepsilon\to 0}\frac{1}{n}\ln\frac{f^n(x_0 + \varepsilon) - f^n(x_0)}{\varepsilon} \\
    &= \lim_{n\to\infty}\frac{1}{n}\ln\frac{df^n(x)}{dx}\bigg|_{x=x_0} \\
    &= \lim_{n\to\infty}\frac{1}{n}\ln\left(\frac{df^n(x)}{df^{n-1}(x)}\bigg|_{f^{n-1}(x)=x_{n-1}}\frac{df^{n-1}(x)}{df^{n-2}(x)}\bigg|_{f^{n-2}(x)=x_{n-2}}\;\hdots\;\frac{df^2(x)}{df(x)}\bigg|_{f(x)=x_1}\frac{df(x)}{dx}\bigg|_{x=x_0}\right) \\
    &= \lim_{n\to\infty}\frac{1}{n}\ln\left(\frac{df(x)}{dx}\bigg|_{x=x_{n-1}}\frac{df(x)}{dx}\bigg|_{x=x_{n-2}}\;\hdots\;\frac{df(x)}{dx}\bigg|_{x=x_1}\frac{df(x)}{dx}\bigg|_{x=x_0}\right) \\
    &= \lim_{n\to\infty}\frac{1}{n}\ln\prod_{i=0}^{n-1}\frac{df(x)}{dx}\bigg|_{x=x_i} \\
    &= \lim_{n\to\infty}\frac{1}{n}\sum_{i=0}^{n-1}\ln\frac{df(x)}{dx}\bigg|_{x=x_i} 
\end{align*}
\begin{equation}
    \lambda = \lim_{n\to\infty}\frac{1}{n}\sum_{i=0}^{n-1}\ln |f'(x_i)| \label{eq:4}
\end{equation}
where the absolute value was implied for cleaner notation. This demonstrates that the Lyapunov exponent is simply the limit of the average of $\ln |f'(x_i)|$ for all steps of the system.

\subsection{Chaotic Behavior of the Logistic Map}

To calculate the Lyapunov exponent of the logistic map numerically, we can add the following R code:
\begin{verbatim}
    epsilon <- 1e-8
    
    x_temp <- numeric(n)
    x_diff <- numeric(n)
    lyapunov_sum <- log(abs((r * (x_0 + epsilon) * 
        (1 - (x_0 + epsilon)) - x[1]) / epsilon))
    
    for (i in 1:(n-1)) {
      x_temp[i+1] = r * (x[i] + epsilon) * (1 - (x[i] + epsilon))
      x_diff[i] <- abs(x_temp[i+1] - x[i+1])
      lyapunov_sum = lyapunov_sum + log(x_diff[i] / epsilon)
    }
    
    lyapunov_exp <- lyapunov_sum / n
    
    print(paste("Lyapunov exponent: ", lyapunov_exp))
\end{verbatim}
This code approximates $\lambda$ well as we use a large $n=10000$ and a small $\varepsilon = 10^{-8}$. In this code, we approximate $f'(x_i)$ from Equation \ref{eq:4} as $\frac{f(x_i+\varepsilon) - f(x_i)}{\varepsilon}=\frac{f(x_i+\varepsilon) - x_{i+1}}{\varepsilon}$. We define two vectors \verb|x_temp| and \verb|x_diff|, which will contain $f(x_i+\varepsilon)$ and $f(x_i+\varepsilon) - f(x_i)$ for all steps, respectively. Next, we define \verb|lyapunov_sum|, which is the running sum $\sum_{i=0}^{n-1}\ln |f'(x_i)|$, setting it initially equal to $\ln \left|\frac{f(x_0+\varepsilon) - x_1}{\varepsilon}\right|$. We then loop through all of the steps, filling in \verb|x_temp| and \verb|x_diff| and using them to add to \verb|lyapunov_sum| each step. Finally, we divide \verb|lyapunov_sum| by $n$ in accordance with Equation \ref{eq:4} and print the value.

\begin{wrapfigure}[20]{r}{0.5\textwidth}
    \vspace{-0.75cm}
    \centering
    \includegraphics[width=0.5\textwidth]{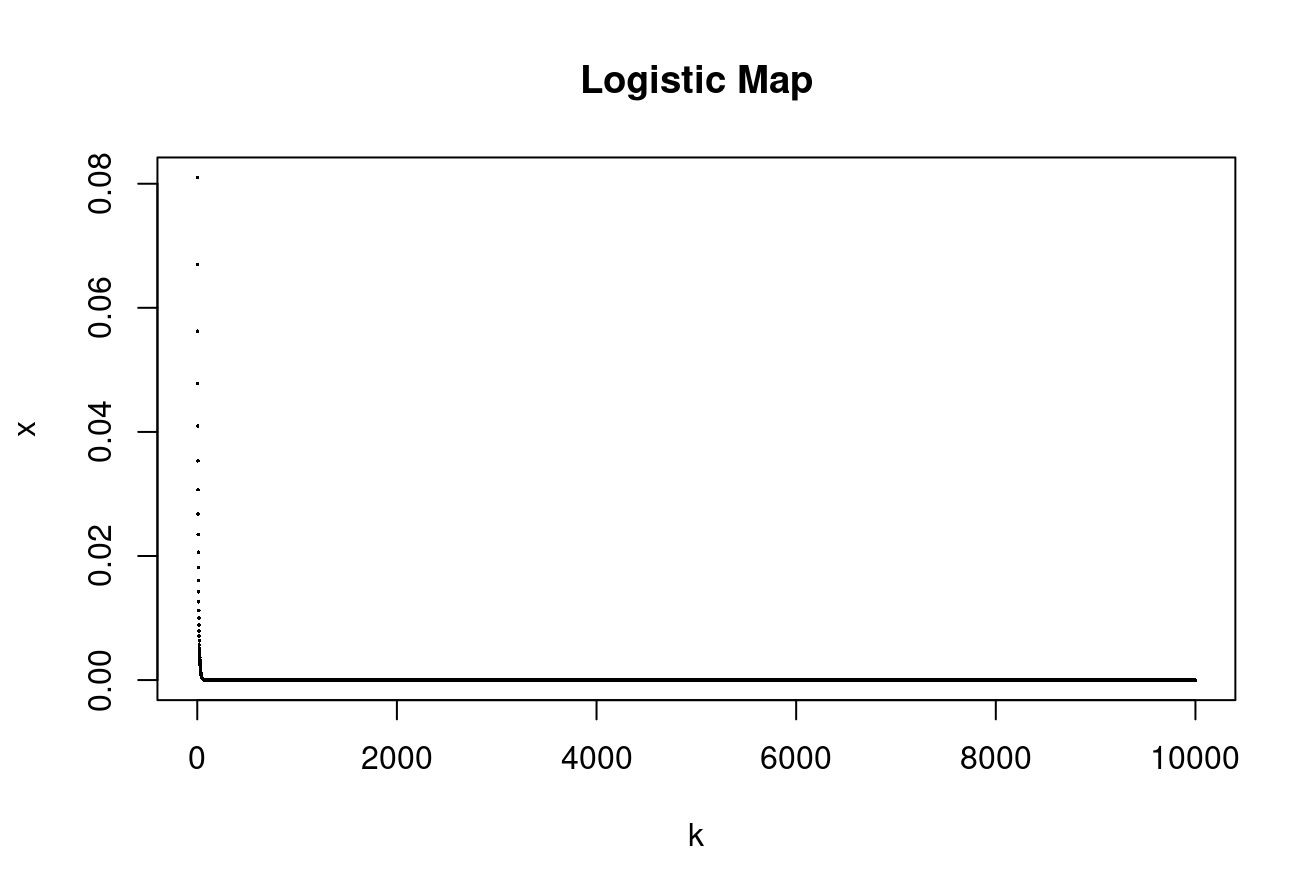}
    \vspace{-0.75cm}
    \caption{Logistic map with $r=0.9$}
    \label{fig:1}
    \includegraphics[width=0.5\textwidth]{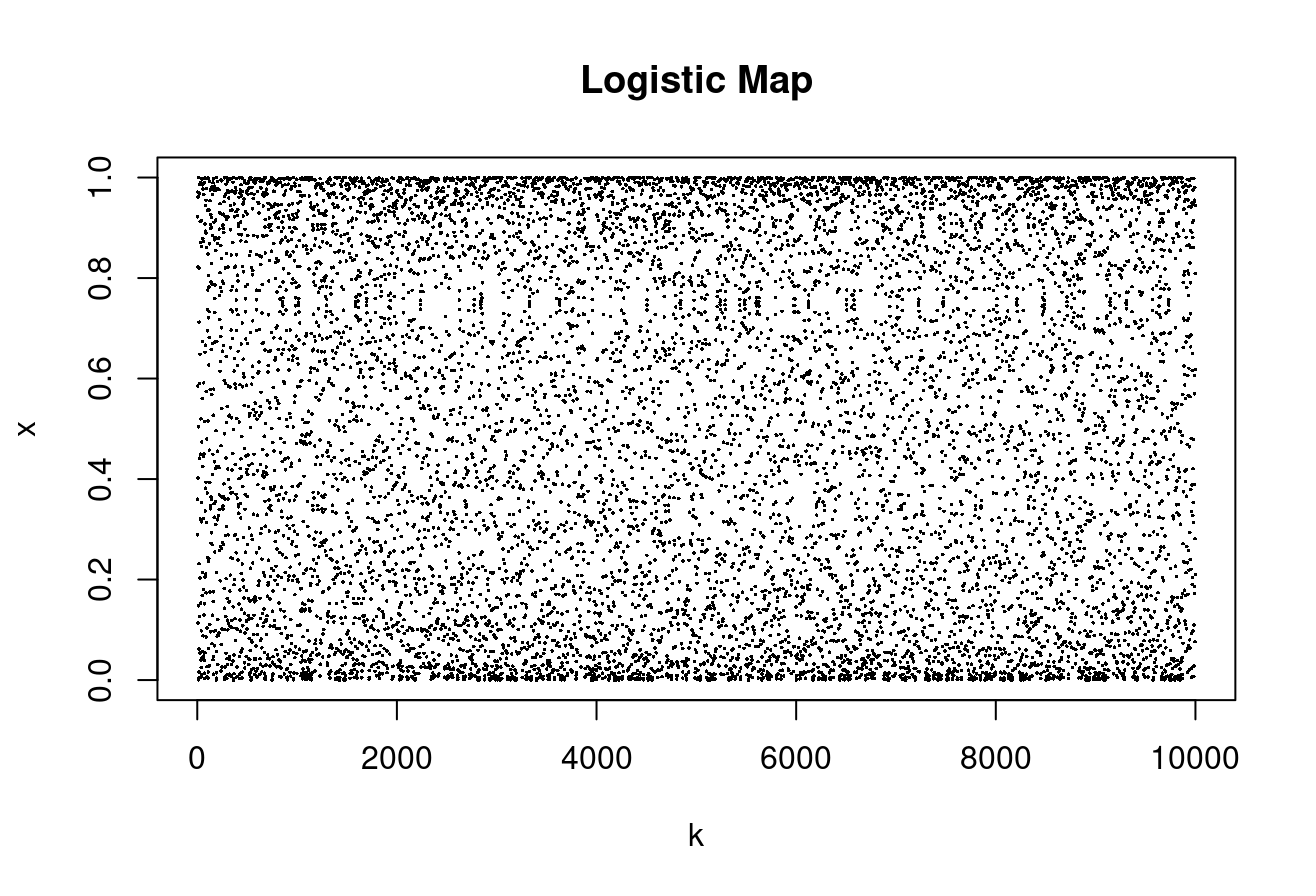}
    \vspace{-0.75cm}
    \caption{Logistic map with $r=4$}
    \label{fig:2}
\end{wrapfigure}

\begin{figure}
    \centering
    \includegraphics[scale = 0.3]{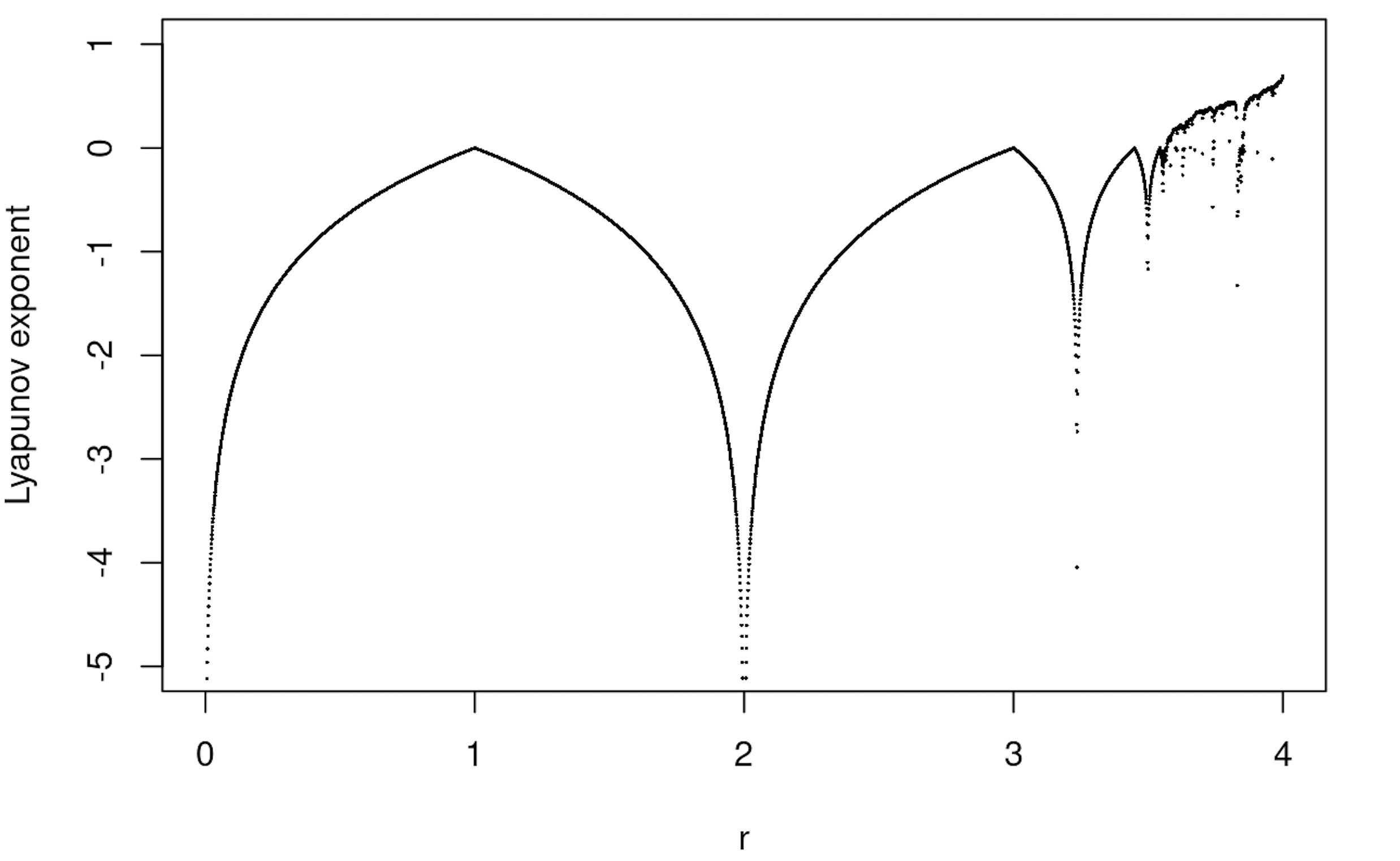}
    \caption{$\lambda$ vs $r$ for $0\leq r\leq4$}
    \label{fig:5}
    \vspace{1cm}
    \includegraphics[scale = 0.3]{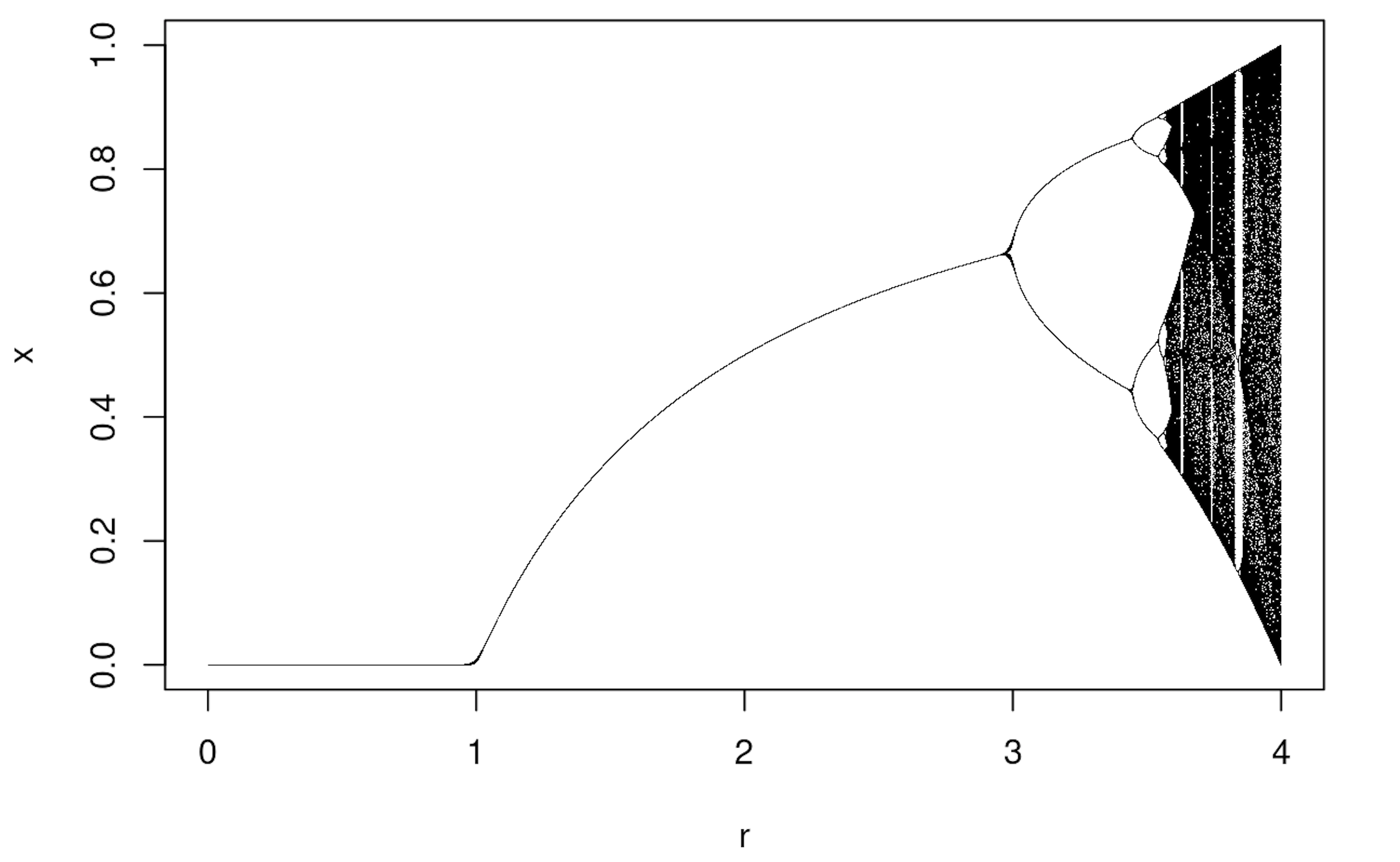}
    \caption{Bifurcation diagram for $0\leq r\leq4$}
    \label{fig:6}
\end{figure}

Let us first consider two values of $r$: 0.9 and 4. For $r=0.9$, we can see in Figure \ref{fig:1} that $x$ begins at its initial value of $x_0=0.1$ then decreases and asymptotically approaches 0. We can assume that if the initial conditions were changed, $x$ would still decrease and asymptotically approach 0, meaning there is little dependence on initial conditions and the system is likely not chaotic. However, for $r=4$, we can see that $x$ is far more complex and appears to be random. We can predict that this system depends on initial conditions because it appears to be chaotic. Using the R code, we find that $\lambda$ for $r=0.9$ is $\lambda = -0.1055$, which indicates that the system is not chaotic, and $\lambda$ for $r=4$ is $\lambda = 0.6932$, which indicates that the system is chaotic. To get a better idea of the chaotic behaviors of the logistic map, we can examine the Lyapunov exponent for many values of $r$ by looping through our established code:

\begin{verbatim}
    start <- 0
    end <- 4
    step <- 0.001
    r <- seq(start, end, by = step)
    r_num <- (end - start + step) / step
    lyapunov_exp <- numeric(r_num)
    n <- 10000
    epsilon <- 1e-8
    k <- numeric(n)
    x <- numeric(n)
    x_temp <- numeric(n)
    x_diff <- numeric(n)
    
    for(j in 1:r_num){
      x_0 <- 0.1
      k[1] <- 1
      x[1] <- r[j] * x_0 * (1 - x_0)
      
      for (i in 1:(n-1)) {
        k[i+1] = i+1
        x[i+1] <- r[j] * x[i] * (1 - x[i])
      }
      
      lyapunov_sum <- log(abs((r * (x_0 + epsilon) * 
        (1 - (x_0 + epsilon)) - x[1]) / epsilon))
      
      for (i in 1:(n-1)) {
        x_temp[i+1] = r[j] * (x[i] + epsilon) * (1 - (x[i] + epsilon))
        x_diff[i] <- abs(x_temp[i+1] - x[i+1])
        lyapunov_sum = lyapunov_sum + log(x_diff[i] / epsilon)
      }
      lyapunov_exp[j] <- lyapunov_sum / n
    }
    
    plot(r, lyapunov_exp, type = "p", pch = 16, 
        xlab = "r", ylab = "Lyapunov exponent", 
        cex = 0.25, ylim = c(-5, 1))
\end{verbatim}
This code produces a plot of many values of $\lambda$ for different values of $r$, which is shown in Figure \ref{fig:5}. To better understand what the graph of $\lambda$ vs $r$ says about the logistic map, we can examine the bifurcation diagram or orbit diagram of the logistic map, which is shown in Figure \ref{fig:6}. The bifurcation diagram plots the attractors, or what states the system tends to after a long time, as a function of $r$ \cite{strogatz}. To do this numerically, we choose a set of $r$ values to consider and loop through these values. In the loop, we iterate through some transients, for example, $x_1$ to $x_{100}$, to reach a steady state of the logistic map. Then, we calculate the values of the points after these transients, for example, $x_{101}$ to $x_{200}$, and plot these for the value of $r$ they correspond to. The following R code accomplishes this and plots Figure \ref{fig:6}:
\begin{verbatim}
    logistic_map <- function(r, x) {
      r * x * (1 - x)
    }
    
    bifurcation_diagram <- function(r_values, x0, num_iterations, num_transient) {
      n <- num_iterations + num_transient
      results <- vector("list", length(r_values))
      for (i in 1:length(r_values)) {
        r <- r_values[i]
        x <- x_0
        iter_result <- numeric(num_iterations)
        
        for (j in 1:num_transient) {
          x <- logistic_map(r, x)
        }
        for (j in 1:num_iterations) {
          x <- logistic_map(r, x)
          iter_result[j] <- x
        }
        results[[i]] <- iter_result
      }
      
      for (i in 1:length(r_values)) {
        r <- r_values[i]
        points(rep(r, num_iterations), results[[i]], 
            pch = ".", col = "black", cex = 0.5)
      }
    }
    
    r_values <- seq(0, 4, by = 0.001)  
    x_0 <- 0.1   
    num_iterations <- 100   
    num_transient <- 100 
    
    plot(1, type = "n", xlim = c(min(r_values), 
        max(r_values)), ylim = c(0, 1),
        xlab = "r", ylab = "x")
    bifurcation_diagram(r_values, x_0, num_iterations, num_transient)
\end{verbatim}

\begin{wrapfigure}[22]{r}{0.4\textwidth}
    \centering
    \includegraphics[width=0.4\textwidth]{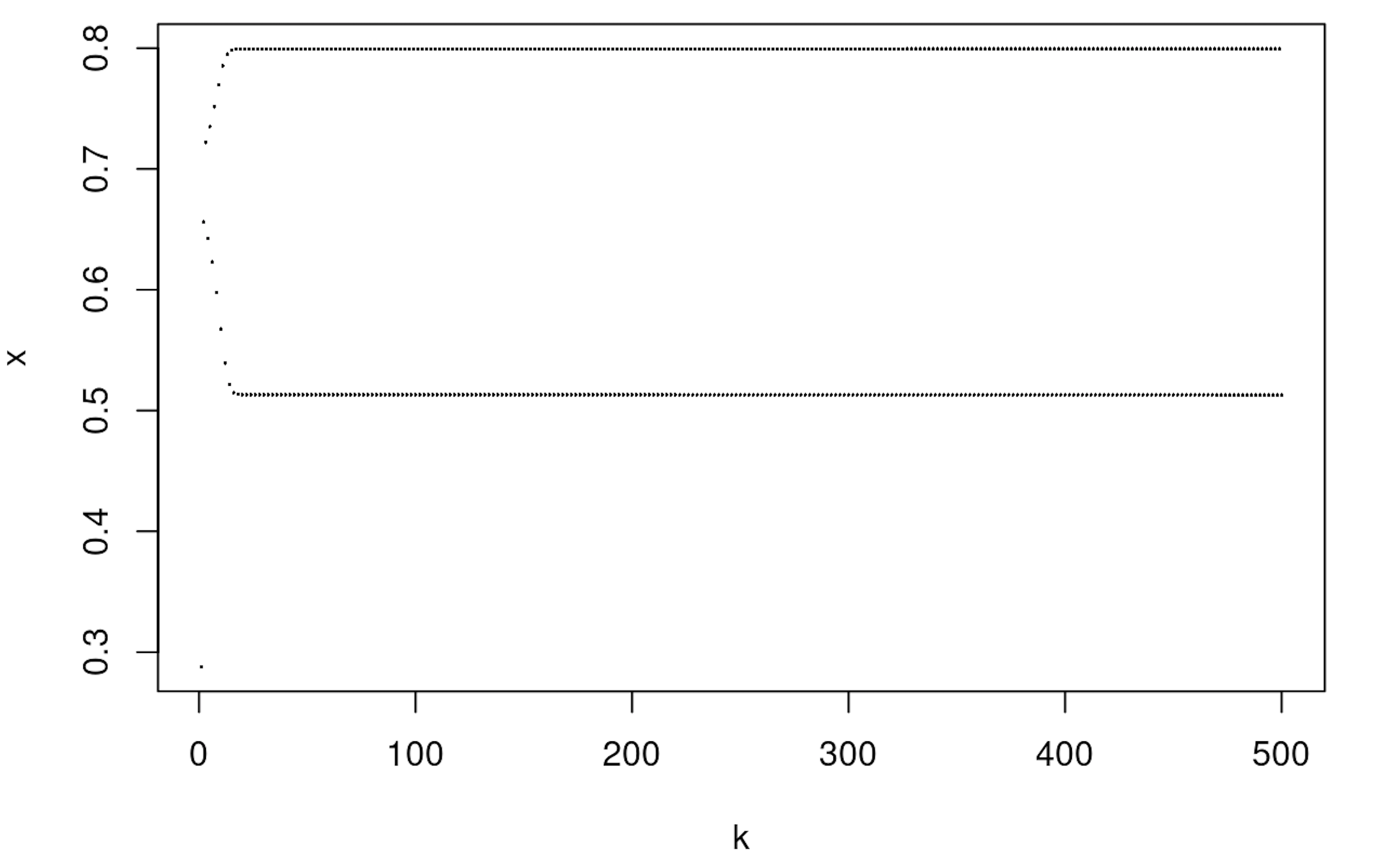}
    \vspace{-0.75cm}
    \caption{Logistic map with $r=3.2$}
    \label{fig:logistic1}
    \vspace{0.25cm}
    \includegraphics[width=0.4\textwidth]{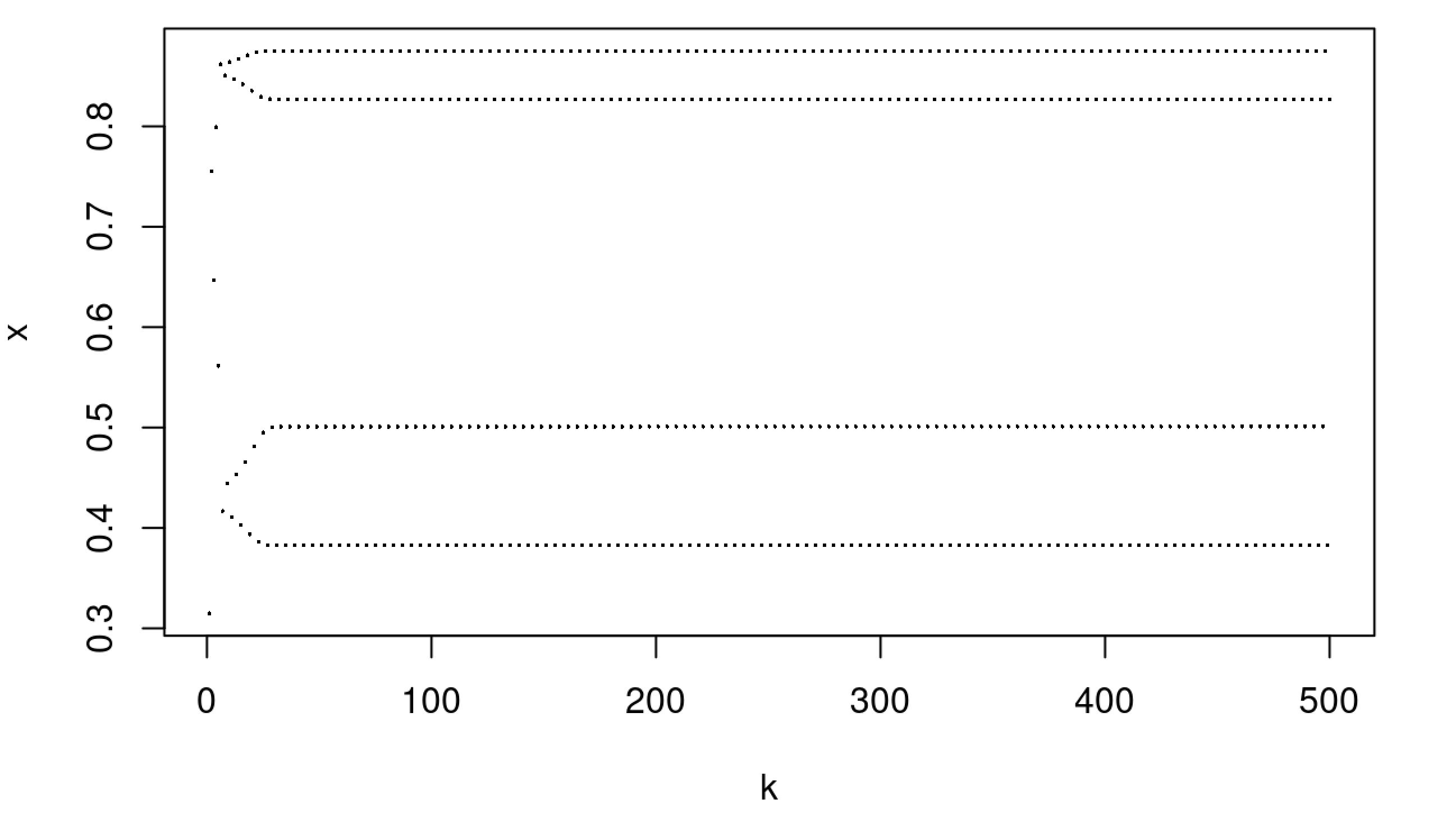}
    \vspace{-0.75cm}
    \caption{Logistic map with $r=3.5$}
    \label{fig:logistic2}
\end{wrapfigure}
\begin{figure}
    \centering
    \includegraphics[scale = 0.3]{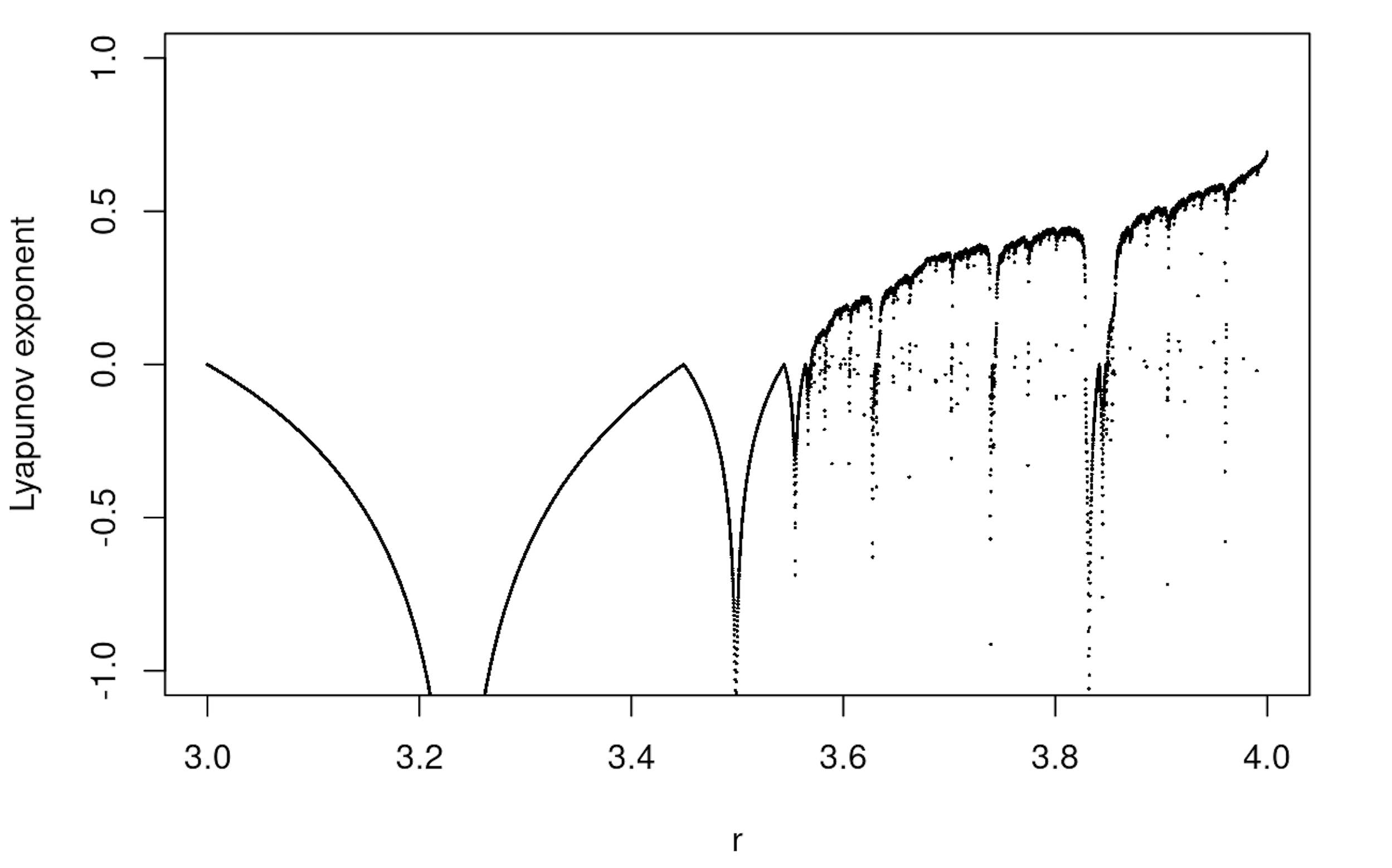}
    \caption{$\lambda$ vs $r$ for $3\leq r\leq4$}
    \label{fig:7}
    \vspace{1cm}
    \includegraphics[scale = 0.35]{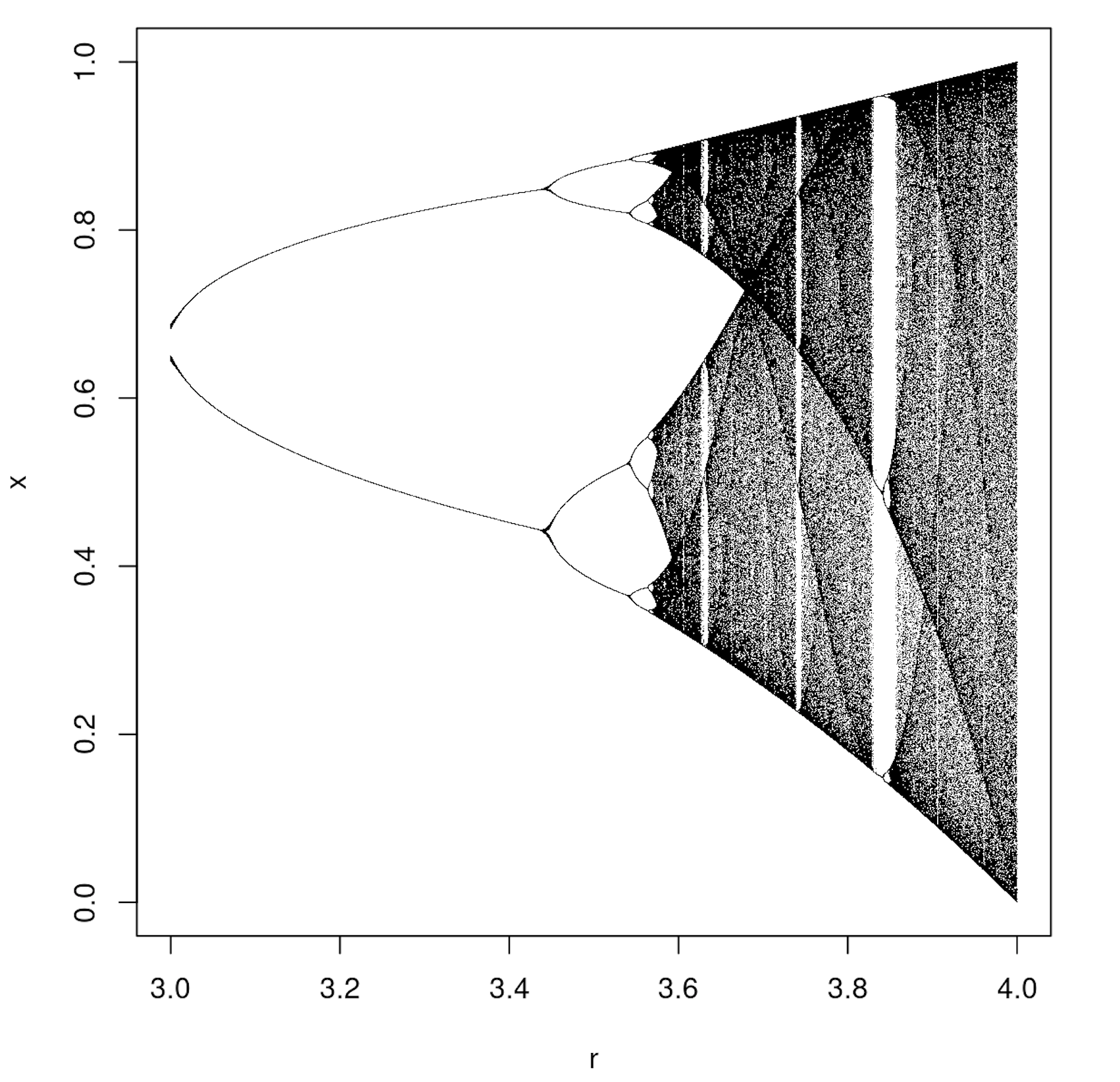}
    \caption{Bifurcation diagram for $3\leq r\leq4$}
    \label{fig:8}
\end{figure}

In Figure \ref{fig:6}, we can see that the logistic map approaches 0 when $r$ is between 0 and 1, which is reflected when $r=0.9$ in Figure \ref{fig:1}. When $r$ is between 0 and 3, the logistic map has an attractor of one point. However, at $r=3$, the branch splits, indicating that the logistic map now has two points as an attractor (Figure \ref{fig:logistic1}). As $r$ increases, the two branches split again to make four branches (Figure \ref{fig:logistic2}), then eight, and so on until chaos is reached and the attractor becomes an infinite set of points (Figure \ref{fig:2}). The point where chaos is reached is $r\approx 3.57$ \cite{strogatz}, which we can see reflected in Figure \ref{fig:7} as the point where $\lambda$ goes above 0 and the point in Figure \ref{fig:8} where the number of attractors blows up. We can also see that the points where the number of attractor points doubles are reflected in the graph of $\lambda$ vs $r$ where $\lambda=0$. In Figure \ref{fig:7}, we can see that past the point of chaos, $\lambda$ becomes discontinuous and there are values of $r$ where the system isn't chaotic, which are called islands of stability. This is easy to see in Figure \ref{fig:8} where there are intervals of $r$ void of many attractor points, most prevalent at $r=1+\sqrt{8}\approx3.83$ \cite{oravec} where there are three attractor points. This island of stability is also clear in Figure \ref{fig:7} by the many points where $\lambda<0$.


\section{Chaos in Multi-Dimensional Systems}

We will now consider chaos for an $n$-dimensional system. Consider a system where every point is defined by a vector $\Vec{r}$ in $n$-dimensional space:
\begin{equation*}
    \Vec{r}=\begin{pmatrix}
        r_1 \\
        r_2 \\
        \vdots \\
        r_n
    \end{pmatrix}
    =\sum_{i=1}^nr_i\Vec{e}_i = r_i\Vec{e}_i
\end{equation*}
where $r_i$ is the component of the $i$th dimension of space $x_i$, $\Vec{e}_i$ is the corresponding basis vector, and we have used the Einstein summation notation. For a discrete-time dynamical system, there is a mapping $(r_1,r_2,\hdots,r_n)\rightarrow(r_1',r_2',\dots,r_n')$ defined by $r_i'=f_i(x_1,x_2,\hdots,x_n)$ for dimension $i$. Given an initial infinitesimal displacement vector $d\Vec{r}=dx_i\Vec{e}_i$,
\begin{align*}
    dx_i' &= \nabla f_i\cdot d\Vec{r} \\
    &= \frac{\partial f_i}{\partial x_1}dx_1 + \frac{\partial f_i}{\partial x_2}dx_2 + \hdots + \frac{\partial f_i}{\partial x_n}dx_n \\
    &= \frac{\partial f_i}{\partial x_j}dx_j
\end{align*}
The Jacobian matrix $J$ is defined as 
\begin{equation*}
    J(\Vec{r}\,) = \begin{pmatrix}
        \frac{\partial f_1}{\partial x_1} & \frac{\partial f_1}{\partial x_2} & \hdots & \frac{\partial f_1}{\partial x_n} \vspace{2.5mm}\\
        \frac{\partial f_2}{\partial x_1} & \frac{\partial f_2}{\partial x_2} & \hdots & \frac{\partial f_2}{\partial x_n} \vspace{2.5mm}\\
        \vdots & \vdots & \ddots & \vdots \vspace{2.5mm}\\
        \frac{\partial f_n}{\partial x_1} & \frac{\partial f_n}{\partial x_2} & \hdots & \frac{\partial f_n}{\partial x_n}
    \end{pmatrix},\hspace{2.5mm}
    J_{ij} = \frac{\partial f_i}{\partial x_j}
\end{equation*}
Substituting,
\begin{equation}
    dx_i'=J_{ij}dx_j \label{eq:5}
\end{equation}
Then, for the vector $d\Vec{r}\,'=dx_i'\Vec{e}_i$,
\begin{equation*}
    \begin{pmatrix}
        dx_1' \\
        dx_2' \\
        \vdots \\
        dx_n'
    \end{pmatrix}=
    \begin{pmatrix}
        J_{1j}dx_j \\
        J_{2j}dx_j \\
        \vdots \\
        J_{nj}dx_j 
    \end{pmatrix}=
    J
    \begin{pmatrix}
        dx_1 \\
        dx_2 \\
        \vdots \\
        dx_n
    \end{pmatrix}
\end{equation*}
\begin{equation}
    d\Vec{r}\,' = J(\Vec{r}\,)d\Vec{r} \label{eq:6}
\end{equation}
Given some initial state $\Vec{r}_0$, we can define the iteration $t$ of the system as $\Vec{r}\,(t)=f^t(\Vec{r}_0)$. Then,
\begin{equation}
    d\Vec{r}\,(t) = J^t(\Vec{r}_0)d\Vec{r}_0 \label{eq:7}
\end{equation}
where, explicitly,
\begin{equation*}
    J^t(\Vec{r}_0) = J(\Vec{r}\,(t-1))\,J(\Vec{r}\,(t-2))\;\hdots\;J(\Vec{r}_0) = J^t
\end{equation*}

Recalling Equation $\ref{eq:3}$, we defined the Lyapunov exponent as $\lambda=\lim_{n\to\infty}\frac{1}{n}\ln\left(\frac{|f^n(x_0 + \varepsilon) - f^n(x_0)|}{\varepsilon}\right)$ for infinitesimal $\varepsilon$ for a one-dimensional discrete-time dynamical system. For our $n$-dimensional case, let us consider an $n$-dimensional initial infinitesimal perturbation vector $d\Vec{r}_0$. After $t$ steps, this vector will evolve to $d\Vec{r}\,(t)$. We can then rewrite Equation \ref{eq:3} for an $n$-dimensional discrete-time dynamical system as
\begin{equation}
    \lambda = \lim_{t\to\infty}\frac{1}{t}\ln\frac{|d\Vec{r}\,(t)|}{|d\Vec{r}_0|} \label{eq:8}
\end{equation}
Substituting Equation \ref{eq:7},
\begin{equation}
    \lambda = \lim_{t\to\infty}\frac{1}{t}\ln\frac{|J^t(\Vec{r}_0)d\Vec{r}_0|}{|d\Vec{r}_0|} = \lim_{t\to\infty}\frac{1}{t}\ln|J^t(\Vec{r}_0)\Vec{u}_0| \label{eq:9}
\end{equation}
where $\Vec{u}_0=d\Vec{r}_0/|d\Vec{r}_0|$ is a unit vector in the direction of $d\Vec{r}_0$. For a vector $\Vec{v} = v_i\Vec{e}_i$,
\begin{equation*}
    \Vec{v}\,^\intercal\Vec{v}=\begin{pmatrix}
        v_1 & v_2 & ... & v_n
    \end{pmatrix}
    \begin{pmatrix}
        v_1 \\
        v_2 \\
        \vdots \\
        v_n
    \end{pmatrix} =
    v_1^2 + v_2^2 + \hdots + v_n^2 = |\Vec{v}|^2
\end{equation*}
Thus, $|J^t(\Vec{r}_0)\Vec{u}_0| = [(J^t\Vec{u}_0)^\intercal J^t\Vec{u}_0]^{1/2} = [\Vec{u}_0^\intercal J^{t\intercal}J^t\Vec{u}_0]^{1/2}$ because $(AB)^\intercal = B^\intercal A\transpose$ \cite[p. 101]{linear}. We can then rewrite Equation \ref{eq:9} as
\begin{equation}
    \lambda = \lim_{t\to\infty}\frac{1}{2t}\ln[\Vec{u}_0^\intercal J^{t\intercal}J^t\Vec{u}_0] = \lim_{t\to\infty}\frac{1}{2t}\ln[\Vec{u}_0^\intercal H_t(\Vec{r}_0)\Vec{u}_0] \label{eq:10}
\end{equation}
where we have defined $H_t(\Vec{r}_0) = J^{t\intercal}(\Vec{r}_0)\,J^t(\Vec{r}_0) = H_t$. Taking $\Vec{u}_0$ to be a normalized eigenvector of $H_t$ $\Vec{w}_i$, we can rewrite Equation \ref{eq:10} as
\begin{equation}
    \lambda_i = \lim_{t\to\infty}\frac{1}{2t}\ln[\Vec{w}_i^\intercal H_t\Vec{w}_i] = \lim_{t\to\infty}\frac{1}{2t}\ln[\Vec{w}_i^\intercal \mu_i\Vec{w}_i] = \lim_{t\to\infty}\frac{1}{2t}\ln \mu_i \label{eq:11}
\end{equation}
where $\mu_i$ is an eigenvalue of $H_t$ and the repeated index $i$ represents $\mu_i$ being matched with its associated eigenvector and is not being summed over. For any matrix $A$, $(A^\intercal A)^\intercal = A^\intercal(A^\intercal)^\intercal = A^\intercal A$, meaning $H_t$ is a real Hermitian or symmetric matrix and has $n$ real eigenvalues (counting multiplicities) \cite[p. 397]{linear}. Therefore, an $n$-dimensional time-discrete dynamical system has $n$ or less Lyapunov exponents. This means that $\lambda$ is a spectrum of Lyapunov exponents $\{\lambda_1,\lambda_2,\hdots,\lambda_n\}$. These Lyapunov exponents and their associated eigenvalues $\mu_i$ are ordered such that $\lambda_1\geq\lambda_2\geq\hdots\geq\lambda_n$ \cite{sano}. Mirroring the one-dimensional system, an $n$-dimensional system is considered chaotic if it contains at least one Lyapunov exponent that is positive, or $\lambda_1>0$, and the magnitude of $\lambda_1$ reflects how quickly a perturbation will be magnified \cite{wolf}. 

Let us now consider starting with a general direction $\Vec{u}_0$. Because $H_t$ is a symmetric matrix, it is orthogonally diagonalizable \cite[p. 398]{linear}, so we can represent $\Vec{u}_0$ as a linear combination of its orthonormal eigenvectors $\Vec{w}_i$
\begin{equation*}
    \Vec{u}_0=a_i\Vec{w}_i
\end{equation*}
where we are using Einstein summation notation again. We can then represent the argument of the natural log in Equation $\ref{eq:10}$ as
\begin{align*}
    \Vec{u}_0^\intercal H_t\Vec{u}_0 &= (a_i\Vec{w}_i)^\intercal H_t\,a_i\Vec{w}_i \\
    &= a_i^2\Vec{w}_i^\intercal H_t\Vec{w}_i \\
    &= a_i^2\mu_i
\end{align*}
Solving Equation \ref{eq:11} for $\mu_i$,
\begin{equation}
    \mu_i \approx e^{2t\Tilde{\lambda}_i} \label{eq:12}
\end{equation}
where $\Tilde{\lambda}_i$ is the approximate Lyapunov exponent for large values of $t$. Substituting,
\begin{align*}
    \Vec{u}_0^\intercal H_t\Vec{u}_0 &\approx a_i^2e^{2t\Tilde{\lambda}_i} \\
    &= a_1^2e^{2t\Tilde{\lambda}_1}+a_2^2e^{2t\Tilde{\lambda}_2}+\hdots+a_n^2e^{2t\Tilde{\lambda}_n} \\
    &\approx a_1^2e^{2t\Tilde{\lambda}_1}
\end{align*}
because the sum is dominated by the first term when $t$ is large since $\lambda_1$ is the largest Lyapunov exponent. Substituting into Equation \ref{eq:10},
\begin{align*}
    \lambda &\approx\lim_{t\to\infty}\frac{1}{2t}\ln(a_1^2e^{2t\Tilde{\lambda}_1}) \\
    &= \lim_{t\to\infty}\left(\frac{\ln a_1^2}{2t} + \frac{\ln{e^{2t\Tilde{\lambda}_1}}}{2t}\right) \\
    &= \Tilde{\lambda}_1
\end{align*}
This means that an arbitrary $\Vec{u}_0$ will almost always produce $\lambda_1$ for large values of $t$. If we want to find $\lambda_2$, we can simply pick a vector where $a_1=0$, which will produce
\begin{align*}
    \Vec{u}_0^\intercal H_t\Vec{u}_0 &\approx a_i^2e^{2t\Tilde{\lambda}_i} \\
    &= a_2^2e^{2t\Tilde{\lambda}_2}+a_3^2e^{2t\Tilde{\lambda}_3}+\hdots+a_n^2e^{2t\Tilde{\lambda}_n} \\
    &\approx a_2^2e^{2t\Tilde{\lambda}_2} \\
\end{align*}
\vspace{-1cm}
\begin{equation*}
    \lambda\approx\lim_{t\to\infty}\frac{1}{2t}\ln(a_2^2e^{2t\Tilde{\lambda}_2}) = \Tilde{\lambda}_2
\end{equation*}
We can then carry on and, in theory, calculate all $n$ Lyapunov exponents. 

\subsection{The H\'enon Map}

As an example of an $n$-dimensional discrete-time chaotic system, let us consider the H\'enon map, a two-dimensional discrete-time dynamical system with mapping 
\begin{align*}
    x' &= f_x(x, y) = 1+y-ax^2 \\
    y' &= f_y(x, y) = bx
\end{align*}
where $a$ and $b$ are parameters \cite{socolar}. The Jacobian $J(\Vec{r}\,)$ of the H\'enon map for $\Vec{r}=x\Vec{e}_x+y\Vec{e}_y$ is 
\begin{equation}
    J(\Vec{r}\,) = \begin{pmatrix}
        \frac{\partial f_x}{\partial x} & \frac{\partial f_x}{\partial y} \vspace{2.5mm}\\
        \frac{\partial f_y}{\partial x} & \frac{\partial f_y}{\partial y}
    \end{pmatrix} = 
    \begin{pmatrix}
        -2x & 1 \\
        b & 0
    \end{pmatrix}
    \label{eq:13}
\end{equation}
If we define $\Vec{r}\,(t) = \Vec{r}_t = x_t\Vec{e}_x + y_t\Vec{e}_y$ for the H\'enon map, we can write Equation \ref{eq:6} as
\begin{equation*}
    \begin{pmatrix}
        dx_t \\
        dy_t
    \end{pmatrix}
    = \begin{pmatrix}
        -2x_{t-1} & 1 \\
        b & 0
    \end{pmatrix}
    \begin{pmatrix}
        dx_{t-1} \\
        dy_{t-1}
    \end{pmatrix}
\end{equation*}
Equation \ref{eq:7} reads
\begin{align*}
    \begin{pmatrix}
        dx_t \\
        dy_t
    \end{pmatrix}
    &= \begin{pmatrix}
        -2x_{t-1} & 1 \\
        b & 0
    \end{pmatrix}
    \begin{pmatrix}
        -2x_{t-2} & 1 \\
        b & 0
    \end{pmatrix}\hdots
    \begin{pmatrix}
        -2x_{0} & 1 \\
        b & 0
    \end{pmatrix}
    \begin{pmatrix}
        dx_{0} \\
        dy_{0}
    \end{pmatrix} \\
    &=
    \prod_{i=0}^{t-1}
    \begin{pmatrix}
        -2x_{i} & 1 \\
        b & 0
    \end{pmatrix}
    \begin{pmatrix}
        dx_{0} \\
        dy_{0}
    \end{pmatrix}
\end{align*}
Notice that 
\begin{equation*}
    \det(J(\Vec{r}_t)) = \begin{vmatrix}
        -2x_{t} & 1 \\
        b & 0
    \end{vmatrix}
    =-b
\end{equation*}
for all $t$. Because $\det(AB) = \det(A)\det(B)$ \cite[p. 175]{linear}, we can then write
\begin{align*}
        \det(J^t(\Vec{r}_0)) &= \det(J(\Vec{r}_{t-1})\,J(\Vec{r}_{t-2})\;\hdots\;J(\Vec{r}_0)) \\
        &= \det(J(\Vec{r}_{t-1}))\,\det(J(\Vec{r}_{t-2}))\;\hdots\;\det(J(\Vec{r}_0))) \\
        &= (-b)(-b)\hdots (-b) \\
        &= (-1)^tb^t
\end{align*}
Since $\det(A\transpose) = \det(A)$ \cite[p. 174]{linear}, we can write
\begin{align*}
    \det(H_t(\Vec{r}_0)) &= \det(J^{t\intercal}(\Vec{r}_0)\,J^t(\Vec{r}_0)) \\
    &= \det(J^{t\intercal}(\Vec{r}_0))\,\det(J^t(\Vec{r}_0)) \\
    &= \det(J^t(\Vec{r}_0))\,\det(J^t(\Vec{r}_0)) \\
    &= (-1)^tb^t\,(-1)^tb^t \\
    &= b^{2t}
\end{align*}
Since $\det(A) = \prod_{i=1}^n\nu_i$ where $\nu_i$ is an eigenvalue of $A$ \cite[p. 282]{linear}, by Equation \ref{eq:12}
\begin{align*}
    \det(H_t(\Vec{r}_0)) &= \mu_1\mu_2  \\
    &\approx e^{2t\Tilde{\lambda}_1}e^{2t\Tilde{\lambda}_2}  \\
    &= e^{2t(\Tilde{\lambda}_1+\Tilde{\lambda}_2)} = b^{2t} 
\end{align*}
\begin{equation}
    \Tilde{\lambda}_1+\Tilde{\lambda}_2\approx\frac{1}{2t}\ln(b^{2t}) = \ln b \label{eq:14}
\end{equation}
meaning a numerical calculation of one Lyapunov exponent will immediately yield the value of the other \cite{ott}.

\begin{figure}[t]
    \centering
    \includegraphics[scale=0.25]{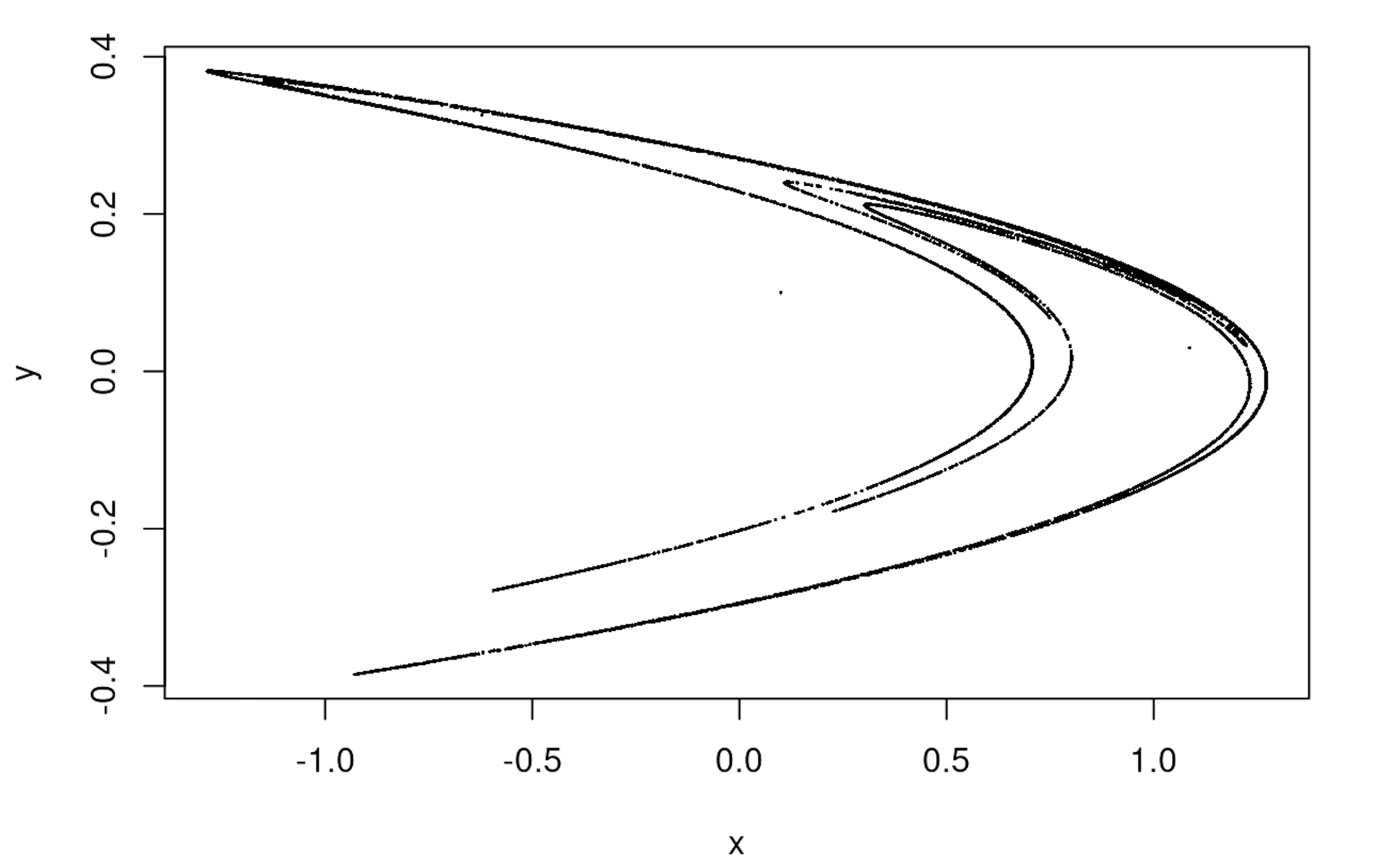}
    \caption{H\'enon map}
    \label{fig:henon}
\end{figure}

Generally, the Lyapunov spectrum cannot be calculated analytically \cite{benettin}, including for the H\'enon map, so we usually resort to numerical techniques. To do this, we must first model the H\'enon map itself. The following code accomplishes this in R:
\begin{verbatim}
    a <- 1.4
    b <- 0.3
    
    t <- 1000000
    
    x_0 <- 0.1
    y_0 <- 0.1
    
    x <- numeric(t)
    y <- numeric(t)
    x[1] <- 1 + y_0 - a * x_0^2
    y[1] <- b * x_0
    
    for (i in 1:(t-1)) {
      x[i+1] <- 1 + y[i] - a * x[i]^2 
      y[i+1] <- b * x[i]
    }

    plot(x, y, type = "p", pch = 16, 
        xlab = "x", ylab = "y", 
        cex = 0.25, main = "Henon Map")
\end{verbatim}
Here, we set the parameters $a$ and $b$ to $1.4$ and $0.3$, respectively, and the number of steps $t$ to 1000000. We then define the set of $x$ and $y$ values as a vector of length $t$ and set the initial point to be $(0.1,\,0.1)$. Finally, we loop through all the steps and apply the H\'enon mapping to fill the $x$ and $y$ vectors and plot the results.

\subsection{A Periodic Renormalization Method}

There are many methods to calculate the Lyapunov spectrum, but one of the most used methods is described by Benettin et al. \cite{benettin}, which we will discuss first. Notice that, although Equation \ref{eq:11} would probably be best for a system that $\lambda$ could be calculated analytically, Equation \ref{eq:9} is best for a numerical calculation because the Jacobian $J(\Vec{r}\,)$ is known, we can simply apply it $t$ times to an arbitrary initial normalized vector $\Vec{u}_0$, producing vectors $\vec{a}_i$. However, for a chaotic system, $\lambda_1$ will be positive, meaning any unit vector will grow to have a very large magnitude after a large number of steps $t$ because any difference in initial conditions is exponentially magnified. Hence, we should periodically renormalize the vector $\vec{a}_{i\tau} = J^\tau\vec{a}_{(i-1)\tau}$ every number of steps $\tau$. In other words, we start with our unit vector $\Vec{u}_0$, apply the Jacobian $\tau$ times, normalize the vector $\Vec{r}_{\tau}$, apply the Jacobian to the normalized vector $\vec{a}_{\tau}/|\vec{a}_{\tau}|$ $\tau$ times, normalize the vector $\vec{a}_{2\tau}$, and so on $s$ times. Then, $t=s\tau$. To see how we can calculate $\Tilde{\lambda}_1$ this way, we have to modify Equation \ref{eq:9}. Setting $\Vec{u}_0=\vec{a}_0$,
\begin{align*}
    |J^t\Vec{u}_0| &= |J^t\vec{a}_0| = \left|\prod_{i=s-1}^{0} J^{\tau}(\Vec{r}_{i\tau})\,\vec{a}_0\right| = \left|\prod_{i=s-1}^{0} J^{\tau}\vec{a}_0\right| \\
    &= \frac{|J^{\tau}\vec{a}_0|}{|J^{\tau}\vec{a}_0|}\left|\prod_{i=s-1}^{0} J^{\tau}\vec{a}_0\right|
    = |J^{\tau}\vec{a}_0|\left|\frac{\prod_{i=s-1}^{0} J^{\tau}\vec{a}_0}{|J^{\tau}\vec{a}_0|}\right|
    = |J^{\tau}\vec{a}_0|\left|\prod_{i=s-1}^{1}J^{\tau}\vec{a}_{\tau}\right| \\
    &= |J^{\tau}\vec{a}_0|\frac{|J^{\tau}\vec{a}_{\tau}|}{|J^{\tau}\vec{a}_{\tau}|}\left|\prod_{i=s-1}^{1}J^{\tau}\vec{a}_{\tau}\right|
    = |J^{\tau}\vec{a}_0||J^{\tau}\vec{a}_{\tau}|\left|\frac{\prod_{i=s-1}^{1} J^{\tau}\vec{a}_{\tau}}{|J^{\tau}\vec{a}_{\tau}|}\right|
    = |J^{\tau}\vec{a}_0||J^{\tau}\vec{a}_{\tau}|\left|\prod_{i=s-1}^{2}J^{\tau}\vec{a}_{2\tau}\right| \\
    &\vdotswithin{=} 
\end{align*}
\begin{align*}
    &= \prod_{i=0}^{s-2}|J^{\tau}\vec{a}_{i\tau}|\frac{|J^\tau\vec{a}_{(s-1)\tau}|}{|J^\tau\vec{a}_{(s-1)\tau}|}|J^{\tau}\vec{a}_{(s-1)\tau}| = \prod_{i=0}^{s-1}|J^{\tau}\vec{a}_{i\tau}|\left|\frac{J^{\tau}\vec{a}_{(s-1)\tau}}{|J^{\tau}\vec{a}_{(s-1)\tau}|}\right| = \prod_{i=0}^{s-1}|J^{\tau}\vec{a}_{i\tau}| \\
    &= \prod_{i=1}^{s}|\vec{a}_{i\tau}|
\end{align*}
Let $\alpha_i = |\vec{a}_{i\tau}|$. Then,
\begin{equation*}
    |J^t\Vec{u}_0| = \prod_{i=1}^{s}\alpha_i
\end{equation*}
Substituting into Equation \ref{eq:9},
\begin{align*}
    \lambda_1 &= \lim_{t\to\infty}\frac{1}{t}\ln|J^t\Vec{u}_0| \\
    &= \lim_{s\to\infty}\frac{1}{s\tau}\ln\prod_{i=1}^{s}\alpha_i \\
    &= \lim_{s\to\infty}\frac{1}{s\tau}\sum_{i=1}^{s}\ln\alpha_i
\end{align*}
For large values of $s$,
\begin{equation}
    \Tilde{\lambda}_1 = \frac{1}{s\tau}\sum_{i=1}^{s}\ln\alpha_i \label{eq:15}
\end{equation}
Notice that $\alpha_i$ is the factor that we use to normalize each $\vec{a}_{i\tau}$. Thus, we can store the magnitude of each $\vec{a}_{i\tau}$ used to normalize it for the subsequent $\tau$ Jacobian iterations and average their natural logs according to Equation \ref{eq:14} at the end.

We can implement this method into R for the H\'enon map by adding the following code:
\begin{verbatim}
    tau <- 50
    s <- t/tau
    
    u_0 <- matrix(c(1, 0), nrow = 2)
    J_0 <- matrix(c(-2*x_0, b, 1, 0), nrow = 2, byrow = FALSE)
    J <- function(i) {
      jacobian <- matrix(c(-2*x[i], b, 1, 0), nrow = 2, byrow = FALSE)
    }
    a <- list()
    alpha <- numeric(s)
    
    a[[1]] <- J_0 %*% u_0
    for(j in 1:(tau-1)){
      a[[j+1]] <- J(j) %*% a[[j]]
    }
    alpha[1] <- sqrt(a[[tau]][1, 1]^2 + a[[tau]][2, 1]^2)
    a[[tau]] <- a[[tau]] / alpha[1]
    lyapunov_sum <- log(alpha[1])
    
    for(i in 2:s){
      for(j in ((i-1)*tau):(i*tau-1)){
        a[[j+1]] <- J(j) %*% a[[j]]
      }
      alpha[i] <- sqrt(a[[i*tau]][1, 1]^2 + a[[i*tau]][2, 1]^2)
      a[[i*tau]] <- a[[i*tau]] / alpha[i]
      lyapunov_sum <- lyapunov_sum + log(alpha[i])
    }
    
    lyapunov_exponent_1 <- lyapunov_sum / t
    lyapunov_exponent_2 <- log(b) - lyapunov_exponent_1
    print(paste(lyapunov_exponent_1, lyapunov_exponent_2))
\end{verbatim}
In this code, we first define a $\tau$ and $s$. Then, we define our initial unit perturbation vector $\Vec{u}_0 = \begin{pmatrix} 1 \\ 0 \end{pmatrix}$, which deviates from the initial conditions $\Vec{r}_0 = \begin{pmatrix} 0.1 \\ 0.1 \end{pmatrix}$, the initial Jacobian and a general Jacobian according to Equation $\ref{eq:13}$. We also make a list \verb|a| representing the set of vectors $\vec{a}_i$ and an array \verb|alpha| representing the set of normalization factors $\alpha_i$. We then define $\Vec{a}_1$ by its definition $J(\Vec{r}_0)\vec{u}_0$ and subsequent vectors $\Vec{a}_{j+1} = J(\Vec{r}_j)\vec{a}_j$ up to $j=\tau-1$. Next, we define $\alpha_1 = |\vec{a}_{\tau}|$ and define \verb|lyapunov_sum| to be equal to $\ln\alpha_1$, which will be the running sums of all $\alpha_i$ in Equation \ref{eq:15}. Then we loop through all values of $t$, applying $\Vec{a}_{j+1} = J(\Vec{r}_j)\vec{a}_j$ and periodically assigning a value $\alpha_i = |\vec{a}_{i\tau}|$ every $\tau$ iterations, normalizing the vector $\vec{a}_{i\tau}$, and adding $\alpha_i$ to \verb|lyapunov_sum|. After the $s$ normalizations and additions, we define $\tilde{\lambda}_1$ by Equation \ref{eq:15}, $\tilde{\lambda}_2$ by Equation \ref{eq:14}, and print. This code yields $\tilde{\lambda}_1 = 0.089537$ and $\tilde{\lambda}_2 = -1.293510$, which indicates chaotic behavior in the H\'enon map.

\subsection{A QR Factorization Method}

Sandri briefly describes another method to calculate the Lyapunov spectrum \cite{sandri} which uses the ability to numerically do a Gram-Schmidt orthonormalization. To utilize this, we once again consider Equation \ref{eq:9}. First, let us define $J(\Vec{r}_i)=J_i$ for easier notation. Then, we can do a QR factorization of $J_0 = Q_1R_1$ where $Q_1$ is a matrix whose columns form an orthonormal basis for the column space of $J_0$ using a Gram-Schmidt orthonormalization and $R_1$ is an upper triangular matrix \cite[p. 359]{linear}. From $k = 1,\,2,\:\hdots\:,\,t-1$, we can define
\begin{equation*}
    Q_{k+1}R_{k+1} = J_kQ_k
\end{equation*}
which is a recursive sequence, first taking the $Q_1$ matrix from the factorization of $J_0$, multiplying by the next Jacobian $J_1$, factorizing again, and so on. Since an orthogonal matrix $U$ has the property $U^{\intercal}U = UU^{\intercal} = I$ \cite[p. 345]{linear}, it follows that
\begin{equation*}
    J_k = Q_{k+1}R_{k+1}Q_k^{\intercal}
\end{equation*}
Applying this to the matrix $J^t$,
\begin{align*}
    J^t &= J_{t-1}J_{t-2}\,\hdots\,J_1J_0 \\
    &= Q_tR_tQ_{t-1}^{\intercal}Q_{t-1}R_{t-1}Q_{t-2}^{\intercal}\,\hdots\,Q_2R_2Q_1^{\intercal}Q_1R_1 \\
    &= Q_tR_tR_{t-1}\,\hdots\,R_2R_1 \\
    &= Q_t\prod_{i=t}^1R_i \\
    &= Q_tR
\end{align*}
where $\prod_{i=t}^1R_i = R$ since the product of upper triangular matrices results in an upper triangular matrix. Explicitly, 
\begin{equation*}
    R = \begin{pmatrix}
        r_{11} & r_{12} & \hdots & r_{1n} \\
        0 & r_{22} & \hdots & r_{2n} \\
        \vdots & \vdots & \ddots & \vdots \\
        0 & 0 & \hdots & r_{nn}
    \end{pmatrix}
\end{equation*}
We can rewrite Equation \ref{eq:9} as
\begin{equation*}
    \lambda = \lim_{t\to\infty}\frac{1}{t}\ln|Q_tR\Vec{u}_0|  
\end{equation*}
Because the eigenvalues of a triangular matrix are on its diagonal \cite[p. 271]{linear}, if we take $\Vec{u}_0$ to be in the direction of an eigenvector of $R$,
\begin{equation*}
    |Q_tR\Vec{u}_0| = |r_{ii}||Q_t\Vec{u}_0|
\end{equation*}
where $i$ isn't being summed over. Because an orthogonal matrix $U$ has the property $|U\Vec{x}| = |\Vec{x}|$ \cite[p. 345]{linear},
\begin{equation*}
    |r_{ii}||Q_t\Vec{u}_0| = |r_{ii}||\Vec{u}_0| = |r_{ii}|
\end{equation*}
Therefore,
\begin{equation*}
    \lambda_i = \lim_{t\to\infty}\frac{1}{t}\ln |r_{ii}|
\end{equation*}
or for large values of $t$,
\begin{equation}
    \Tilde{\lambda}_i = \frac{1}{t}\ln |r_{ii}| \label{eq:16}
\end{equation}
This gives us a way to calculate the entire Lyapunov spectrum using QR factorization.

If we consider a computational approach to implementing this, for large values of $t$, the values $r_{ii}$ will grow very large. Thus, we cannot find the Lyapunov spectrum by calculating $R$ through direct matrix multiplication. If we say an element of the matrix $R_j$ is $r_{j_{rc}}$, consider multiplying two matrices $R_2$ and $R_1$:
\begin{equation*}
    R_2R_1 = \begin{pmatrix}
        r_{2_{11}} & r_{2_{12}} & \hdots & r_{2_{1n}} \\
        0 & r_{2_{22}} & \hdots & r_{2_{2n}} \\
        \vdots & \vdots & \ddots & \vdots \\
        0 & 0 & \hdots & r_{2_{nn}}
    \end{pmatrix}
    \begin{pmatrix}
        r_{1_{11}} & r_{1_{12}} & \hdots & r_{1_{1n}} \\
        0 & r_{1_{22}} & \hdots & r_{1_{2n}} \\
        \vdots & \vdots & \ddots & \vdots \\
        0 & 0 & \hdots & r_{1_{nn}}
    \end{pmatrix}
    =
    \begin{pmatrix}
        r_{2_{11}}r_{1_{11}} & x_{12} & \hdots & x_{1n} \\
        0 & r_{2_{22}}r_{1_{22}} & \hdots & x_{2n} \\
        \vdots & \vdots & \ddots & \vdots \\
        0 & 0 & \hdots & r_{2_{nn}}r_{1_{nn}}
    \end{pmatrix}
\end{equation*}
The off-diagonal elements are irrelevant to calculating the Lyapunov exponents. Continuing this pattern, 
\begin{equation*}
    \ln |r_{ii}| = \ln\prod_{j=1}^t|r_{j_{ii}}| = \sum_{j=1}^t\ln |r_{j_{ii}}|
\end{equation*}
Substituting into Equation \ref{eq:16},
\begin{equation}
        \Tilde{\lambda}_i = \frac{1}{t}\sum_{j=1}^t\ln |r_{j_{ii}}| \label{eq:17}
\end{equation}

We can now implement this method for the H\'enon map with the following R code: 
\begin{verbatim}
    J_0 <- matrix(c(-2*x_0, b, 1, 0), nrow = 2, byrow = FALSE)
    J <- function(i) {
      jacobian <- matrix(c(-2*x[i], b, 1, 0), nrow = 2, byrow = FALSE)
    }
    QR <- list()
    Q <- list()
    R <- list()
    
    QR[[1]] <- qr(J_0)
    Q[[1]] <- qr.Q(QR[[1]])
    R[[1]] <- qr.R(QR[[1]])
    
    for (i in 1:(t-1)) {
      QR[[i+1]] <- qr(J(i) %*% Q[[i]])
      Q[[i+1]] <- qr.Q(QR[[i+1]])
      R[[i+1]] <- qr.R(QR[[i+1]])
    }
    
    lyapunov_1sum <- 0
    lyapunov_2sum <- 0
    
    for (j in 1:t) {
      lyapunov_1sum <- lyapunov_1sum + log(abs(R[[j]][1,1]))
      lyapunov_2sum <- lyapunov_2sum + log(abs(R[[j]][2,2]))
    }
    
    lyapunov_exponent_1 <- lyapunov_1sum / t
    lyapunov_exponent_2 <- lyapunov_2sum / t
    
    print(paste(lyapunov_exponent_1, lyapunov_exponent_2))
\end{verbatim}
In this code, we define an initial Jacobian matrix $J_0$ and a function for the general Jacobian of the H\'enon map. We also make lists that will contain the matrices $QR$, $Q$, and $R$. We then do a QR factorization according to the definition $J_0 = Q_1R_1$ and store these in the list. Next, we loop through all of the values of $t$ and fill the lists of matrices $Q$ and $R$. Finally, we use Equation \ref{eq:17} to calculate the Lyapunov spectrum with a running sum. Using $t = 100000$, this code yields $\Tilde{\lambda}_1 = 0.089576$ and $\Tilde{\lambda}_2 = -1.293549$. These results are very similar to those calculated by the method of renormalizing the perturbation vector, but they required fewer steps, and we were able to calculate both exponents directly.

\newpage
\bibliographystyle{abbrv}
\bibliography{refs}

\end{document}